\documentclass[11pt]{amsart}
\usepackage[margin=1in]{geometry}
\usepackage{amsmath,amssymb,amsthm,mathtools}
\usepackage{microtype,needspace}
\usepackage[colorlinks=true,linkcolor=blue,citecolor=blue,urlcolor=blue]{hyperref}
\newtheorem{theorem}{Theorem}[section]
\newtheorem{lemma}[theorem]{Lemma}
\newtheorem{proposition}[theorem]{Proposition}
\theoremstyle{remark}

\newcommand{\Ric}{\operatorname{Ric}}
\newcommand{\tr}{\operatorname{tr}}
\newcommand{\diam}{\operatorname{diam}}
\newcommand{\osc}{\operatorname{osc}}
\newcommand{\Vol}{\operatorname{Vol}}
\newcommand{\ddc}{i\partial\bar\partial}
\newcommand{\N}{\mathcal N}
\newcommand{\R}{\mathbb R}
\title[The Calabi--Yau theorem by K\"ahler--Ricci flow]{A geometric K\"ahler--Ricci flow proof\\
of the Calabi--Yau theorem}
\author{Bennett Chow and Yuan Liao}
\address{Department of Mathematics, University of California San Diego,
La Jolla, CA 92093, USA}
\date{November 8, 2026}
\begin{document}
\begin{abstract}
We give a geometric proof that every K\"ahler class on a compact K\"ahler
manifold with vanishing real first Chern class contains a unique Ricci-flat
metric, and that the K\"ahler--Ricci flow converges smoothly to this metric.
The fixed K\"ahler class and the bounded logarithmic volume density control
the mean pairwise distance. Chow's scalar-curvature decay makes available,
at a fixed spatial scale and all sufficiently late times, the relative-volume
comparison of Jian--Song, derived from Bamler's pointed Nash entropy
estimates. This comparison upgrades the mean-distance bound to a diameter
bound without an absolute lower bound for entropy. A restarted heat-gradient
estimate then gives exponential decay of the Ricci potential. A parabolic
Schwarz estimate and the interior estimates of Sherman--Weinkove yield
smooth exponential convergence. The proof does not invoke the
Calabi--Yau existence theorem or global estimates for a scalar complex
Monge--Amp\`ere equation.
\end{abstract}
\maketitle

\section{Introduction}

Calabi's conjecture \cite{Calabi} asks for K\"ahler metrics with prescribed
Ricci form. Yau's solution \cite{YauCalabi} implies, in particular, that a
compact K\"ahler manifold whose real first Chern class vanishes admits a
unique Ricci-flat K\"ahler metric in every K\"ahler class. Cao
\cite{Cao} subsequently proved this result by K\"ahler--Ricci flow and
established convergence of the flow. His work is the fundamental antecedent
of the argument presented here.

The purpose of this paper is to give another proof of this classical
theorem, with the long-time estimates organized around distance, volume,
and the heat equation for the Ricci potential. The conclusion is the
following familiar form of the Calabi--Yau theorem together with Cao's
convergence theorem.

\begin{theorem}[Calabi--Yau; Cao]\label{thm:main}
Let \((X,J,\omega_0)\) be a connected compact K\"ahler manifold of complex
dimension \(n\geq1\), with
\(c_1(X)=0\) in \(H^2(X;\R)\). Then \([\omega_0]\) contains a unique
Ricci-flat K\"ahler form \(\omega_\infty\). The K\"ahler--Ricci flow
\[
 \partial_t\omega=-2\rho(\omega),\qquad \omega(0)=\omega_0,
\]
exists for all \(t\geq0\). For every integer \(k\geq0\), there are constants
\(C_k,c_k>0\) such that
\[
 \|\omega(t)-\omega_\infty\|_{C^k(g_0)}\leq C_k e^{-c_kt}
 \qquad(t\geq0).
\]
The complex structure and the underlying manifold remain fixed throughout.
\end{theorem}

All manifolds in this paper are nonempty and without boundary. The
hypothesis concerns the real first Chern class; no trivialization of the
canonical bundle is assumed.

\subsection{The argument}
Fix a positive volume form \(\Omega\) with zero Ricci form and total
volume \(\int_XdV_{g_0}\). The appendix constructs it by a linear Poisson
equation. Along the tensor flow,
\(f=\log(\Omega/dV_{g(t)})\) satisfies the heat equation. Its uniform bound
gives two-sided volume-density bounds, and a parabolic Schwarz estimate
then gives existence on every finite time interval. It remains to obtain
estimates that do not deteriorate as time tends to infinity.

The proof then has three steps. First, cohomology bounds
\(\int_X\tr_{g_0}g(t)\,dV_{g_0}\). Applying a segment inequality on the
fixed background \(g_0\) bounds the mean pairwise distance for \(g(t)\).
Second, Chow's estimate \cite{Chow} gives \(|R|\leq C/(1+t)\).
The relative-volume comparison of Jian--Song \cite[Proposition 3.1]{JS},
obtained from Bamler's entropy estimates \cite{Bamler}, then gives a
uniform diameter bound: a distant ball would carry a definite fraction
of the total volume, contradicting the mean-distance bound. Third, a
heat-gradient estimate restarted at arbitrary times contracts
\(\osc f\) by a fixed factor over a fixed time interval. The resulting
exponential decay bounds the time integral of \(f\). The same Schwarz
estimate, now with uniform constants, and the tensor interior estimates
of Sherman--Weinkove \cite{SW} finish the convergence argument.

The order of these steps matters. The diameter estimate precedes the
uniform metric-equivalence estimate. It uses only scalar-curvature decay,
the bounded volume density, and a cohomological integral on the fixed
background. It requires neither a lower Ricci-curvature bound for the
evolving metric nor a time-uniform noncollapsing assumption.

\subsection{Relation to earlier work}
The individual analytic tools have substantial antecedents. We use
Hamilton's tensor Ricci flow \cite{Hamilton}, the gradient method of Chow
\cite{Chow}, the holomorphic-map Bochner calculation underlying Yau's
Schwarz lemma \cite{YauSchwarz}, and the local derivative estimates of
Sherman--Weinkove \cite{SW}, whose proof incorporates a local parabolic
version of Calabi's third-order estimate. The fixed-background distance
argument is an application of the segment-inequality method of Cheeger--Colding
\cite{CC}; we include the short version needed here. 
Bamler's pointed Nash
entropy estimates \cite{Bamler} build on the entropy approach to Ricci flow
initiated by Perelman \cite{Perelman}. 
We use specifically the cancellation
of the entropy factors in a relative-volume comparison. This cancellation
and the resulting comparison are already present in Jian--Song
\cite[Proposition 3.1]{JS}; we reproduce the argument to make the hypotheses
explicit. Our passage from volume concentration to diameter also follows
their strategy \cite[Lemma 4.4 and Proposition 4.1]{JS}, with a first-moment
bound supplying the concentration in the present setting.

Guo--Phong--Song--Sturm \cite[Theorem 1.1]{GPSS} prove diameter bounds for much more
general families of K\"ahler metrics. Their theorem also applies to the
fixed-class family with bounded volume density considered here. Their
proof uses auxiliary complex Monge--Amp\`ere equations and Yau's existence
theorem, so we do not import it into a proof of Ricci-flat existence.
Our diameter argument is instead a short application of the preceding
Riemannian relative-volume comparison to the mean-distance bound.

A detailed exposition of Perelman's uniform scalar-curvature and
diameter estimates for the normalized K\"ahler--Ricci flow on Fano
manifolds is given by \v{S}e\v{s}um and Tian \cite{ST}.
The work of Jian--Song--Tian \cite{JST} gives a further important
comparison: it uses heat-kernel and Harnack methods to reprove Perelman's
scalar-curvature and diameter estimates on Fano manifolds. Their argument
controls the Ricci potential through a Harnack estimate. In the
\(c_1=0\) setting here, we first bound diameter and then use it to prove
decay of the Ricci potential, which solves the heat equation.

For a geometric approach to the surjectivity of the period map for
hyperk\"ahler metrics on K3 surfaces, see Liu \cite{Liu}.

\subsection{Conventions and analytic inputs}
We use real Ricci-flow time throughout:
\begin{equation}\label{eq:flow}
 \partial_tg=-2\Ric(g),\qquad \partial_t\omega=-2\rho(\omega).
\end{equation}
Write \(m=2n\). In holomorphic coordinates,
\[
 \omega=i g_{i\bar j}\,dz^i\wedge d\bar z^j,\qquad
 g=2\operatorname{Re}(g_{i\bar j}\,dz^i\otimes d\bar z^j),\qquad
 \rho=-i\partial\bar\partial\log\det(g_{i\bar j}).
\]
The scalar curvature \(R\), Laplacian \(\Delta\), gradient, Hessian, and
norms are real Riemannian ones. In particular,
\[
 \Delta=2g^{i\bar j}\partial_i\partial_{\bar j},\qquad
 R=2\tr_\omega\rho,\qquad dV_g=\omega^n/n!.
\]
We write \(L=\partial_t-\Delta_{g(t)}\) and
\(\ddc=i\partial\bar\partial\). An unmarked trace between K\"ahler forms
is a complex trace; \(\tr_h g\) for real metrics is a real trace.
All \(C^k(g_0)\) norms use the fixed background metric.

We use tensor Ricci-flow short-time existence, uniqueness, and preservation
of the K\"ahler condition with the original complex structure
\cite{Hamilton,Cao}, together with scalar maximum principles.
Sherman--Weinkove's local regularity theorem and Bamler's three estimates
are stated at their points of use. The remaining calculations are included
below.

\section{The Ricci potential and long-time existence}\label{sec:existence}

\subsection{A fixed volume form and an evolving heat solution}

The compact K\"ahler \(\partial\bar\partial\)-lemma, in the form proved in
Appendix~\ref{app:potential}, gives a smooth real \(f_0\) with
\(\rho(\omega_0)=\ddc f_0\). Adjust its additive constant so that
\[
 \Omega=e^{f_0}dV_{g_0},\qquad \int_X\Omega=\int_XdV_{g_0}=:V.
\]
This is a positive smooth volume form satisfying
\(\rho(\Omega)=0\), where the Ricci form of a volume form is computed
from its density in holomorphic coordinates. It is not the volume form of an assumed
Ricci-flat metric. Along the short-time tensor solution define
\[
 f(x,t)=\log\frac{\Omega}{dV_{g(t)}}.
\]
The volume evolution \(\partial_t dV_g=-R\,dV_g\) and the Ricci-form
identity give
\begin{equation}\label{eq:potential}
 \rho(\omega)=\ddc f,\qquad \Delta f=R,\qquad f_t=R,\qquad Lf=0.
\end{equation}
Consequently, with \(F=\|f_0\|_\infty\),
\begin{equation}\label{eq:density}
 |f|\leq F,\qquad e^{-F}\Omega\leq dV_{g(t)}\leq e^F\Omega.
\end{equation}
Since \([\rho]=2\pi c_1(X)=0\), the class \([\omega(t)]\) equals
\([\omega_0]\), so the volume is \(V\). This fixes the additive
normalization of \(f\) at every time.
In particular
\begin{equation}\label{eq:normalization}
 \int_X e^{-f}\Omega=V,\qquad
 \min_X f\leq0\leq\max_X f,\qquad
 \|f\|_\infty\leq\osc_Xf.
\end{equation}

\subsection{A trace estimate used twice}
We record a tensor estimate which will be used both at finite and infinite
time. Put
\[
 \Psi(x,t)=2\int_0^t f(x,s)\,ds,\qquad S=\tr_{\omega(t)}\omega_0.
\]
Integrating the tensor equation gives
\begin{equation}\label{eq:primitive}
 \ddc\Psi=\omega_0-\omega(t),\qquad
 L\Psi=2(f+n-S).
\end{equation}
The holomorphic identity map
\((X,\omega(t))\longrightarrow(X,\omega_0)\) satisfies the parabolic
Schwarz estimate, obtained from the holomorphic-map Bochner formula
underlying \cite{YauSchwarz},
\begin{equation}\label{eq:schwarz}
 L\log S\leq 2BS
\end{equation}
for a fixed \(B\geq0\), an upper bound for the holomorphic bisectional
curvature of \(\omega_0\). For clarity about the normalization, write \(\tau=2t\) and
\(\Delta_{\mathbb C}=\Delta/2\). For the flow parametrized by \(\tau\),
the Bochner formula reads
\[
 (\partial_\tau-\Delta_{\mathbb C})S
 =-\mathcal E+
   g^{i\bar j}g^{k\bar\ell}(R_0)_{i\bar j k\bar\ell}
 \leq-\mathcal E+BS^2.
\]
Here \(\mathcal E\) is the squared norm of the covariant derivative of the
holomorphic identity map, using the usual complex Hermitian contractions.
The domain Ricci term cancels the inverse-metric time derivative, and
the corresponding Hermitian norm satisfies
\(|\partial S|_{\mathbb C}^2\leq S\mathcal E\).
Taking the logarithm and using \(L=2(\partial_\tau-\Delta_{\mathbb C})\) gives
\eqref{eq:schwarz}.

Choose \(A>B\). Equations \eqref{eq:primitive}--\eqref{eq:schwarz} imply
\begin{equation}\label{eq:Q}
 L(\log S+A\Psi)
 \leq-2(A-B)S+2A(n+F).
\end{equation}
If \(|\Psi|\leq K\) on a time interval starting at zero, the maximum principle
bounds \(S\) by a constant depending only on \(K,A,B,F,n\): at a positive-time
maximum of \(\log S+A\Psi\), one has
\(S\leq\frac{A(n+F)}{A-B}\), while at time zero \(S=n,\Psi=0\).
The bound on \(\Psi\) transfers this maximum bound to \(S\) everywhere:
\[
 S\leq e^{2AK}\max\left\{n,\frac{A(n+F)}{A-B}\right\}.
\]

An upper bound on \(S=\sum_j\lambda_j^{-1}\), together with
\eqref{eq:density}, bounds all the eigenvalues \(\lambda_j\) of \(\omega(t)\)
relative to \(\omega_0\) above and below. Indeed, it bounds each
\(\lambda_j\) below, and the volume ratio bounds \(\prod_j\lambda_j\) above.

\subsection{Continuation on every finite interval}
We use the following precise tensor regularity input:
Sherman--Weinkove \cite[Theorem 1.1 and Corollary 1.2]{SW}, pp.~491--492,
prove that local uniform metric equivalence along K\"ahler--Ricci flow gives
interior bounds for curvature, its derivatives, and all fixed-background
derivatives of the metric at positive elapsed time. Their proof uses a
local parabolic version of Calabi's third-order estimate to control the
connection difference, equivalently the first derivatives of the metric.
Here ``third-order'' refers to derivatives of a local K\"ahler potential.
Curvature and its higher derivatives are then bounded by maximum-principle
arguments, and linear Poisson estimates give higher coordinate derivatives
of the metric. This proof does not invoke Evans--Krylov theory. On slabs
of a fixed positive length, with fixed coordinate balls and a fixed
equivalence constant, these
bounds are uniform and independent of the metric derivatives at the
initial endpoint. The factor two in \eqref{eq:flow} is removed by a constant
time change.

Suppose the maximal smooth existence time \(T\) were finite.
The bound \(|\Psi|\leq2FT\) and \eqref{eq:Q} yield
\(N_T^{-1}\omega_0\leq\omega(t)\leq N_T\omega_0\) for \(t<T\).
Apply the preceding interior estimates on fixed coordinate balls and time
slabs of length \(T/4\) ending at \(t\in[T/2,T)\).
All spatial derivatives are uniformly bounded there. The tensor flow
equation bounds time derivatives in every fixed-background spatial norm.
Thus \(g(t)\) converges smoothly as \(t\uparrow T\) to a positive K\"ahler
metric, from which short-time existence extends the flow. This contradicts
maximality. The flow therefore exists for all \(t\geq0\).

\section{Scalar-curvature decay}\label{sec:scalar}

We give the maximum-principle calculation behind the scalar-decay estimate
of Chow \cite{Chow}, with all constants in the real conventions above. Let \(q=|\nabla f|^2\) and \(H=\nabla^2f\).
The \(J\)-invariant part
\[
 H^+(u,v)=\tfrac12\bigl(H(u,v)+H(Ju,Jv)\bigr)
\]
equals \(\Ric\), by \(\rho=\ddc f\). The two \(J\)-types are orthogonal, so
\(|H|^2\geq|\Ric|^2\).
Bochner's formula, the real Ricci-flow equation, and \eqref{eq:potential}
give
\begin{equation}\label{eq:bochner}
 Lq=-2|H|^2,\qquad L(f^2)=-2q,\qquad LR=2|\Ric|^2.
\end{equation}
Hence
\[
 L(tq+f^2/2)=-2t|H|^2\leq0,\qquad tq\leq F^2/2.
\]
For the scalar upper bound, introduce the nonsingular quantity
\[
 Z=t^2(R+2q)+2tf^2.
\]
Equation \eqref{eq:bochner} gives
\begin{align*}
 LZ&=2tR+2f^2+t^2(2|\Ric|^2-4|H|^2)\\
 &\leq 2tR+2F^2-\frac2m(tR)^2
 \leq \frac m2+2F^2.
\end{align*}
Since \(Z(0)=0\), the maximum principle yields
\(Z\leq(m/2+2F^2)t\). The other terms in \(Z\) are nonnegative,
so \(R\leq(m/2+2F^2)/t\) for \(t>0\).
The inequality \(LR\geq2R^2/m\) and the scalar maximum principle give
\(R\geq-m/(2t)\) for \(t>0\).
Combining these estimates, and using smoothness for \(0\leq t\leq1\), yields
\begin{equation}\label{eq:Chow}
 |R(x,t)|+|\nabla f(x,t)|^2\leq \frac{C}{1+t}.
\end{equation}

\section{The fixed class controls mean distance}\label{sec:mean}

The next lemma uses cohomology only through a fixed-background trace
integral. No curvature bound for \(g(t)\) enters its proof.

\begin{lemma}\label{lem:mean}
Along this flow, for \(\mu_t=V^{-1}dV_{g(t)}\) there is a finite constant
\(L_0>0\), independent of \(t\), such that
\[
 \iint_{X\times X}d_{g(t)}(x,y)\,d\mu_t(x)\,d\mu_t(y)\leq L_0.
\]
\end{lemma}
\begin{proof}
Write \(h=g_0\). The real trace is twice the complex trace; cohomology gives
\begin{equation}\label{eq:trace-integral}
 \int_X\tr_hg(t)\,dV_h
 =\frac{2n}{n!}\int_X\omega(t)\wedge\omega_0^{n-1}=mV.
\end{equation}
For almost every pair \(x,y\), let \(\gamma_{xy}\) be its unique minimizing
\(h\)-geodesic. Since \(|v|_{g(t)}\leq\sqrt{\tr_hg(t)}\,|v|_h\),
\[
 d_{g(t)}(x,y)\leq
 \int_{\gamma_{xy}}\sqrt{\tr_hg(t)}\,ds_h.
\]
The segment-inequality argument of Cheeger--Colding
\cite[Theorem 2.11, p.~199]{CC}, on the
single fixed compact background, gives
\begin{equation}\label{eq:segment}
 \iint\int_{\gamma_{xy}}G\,ds_h\,dV_h(x)dV_h(y)
 \leq C_h\int_XG\,dV_h
\end{equation}
for every nonnegative integrable \(G\).

Here is a direct justification of \eqref{eq:segment}. Split each segment
into two equal halves and use symmetry to retain the half farther from \(x\).
In \(h\)-polar coordinates about \(x\), write the endpoint distance as
\(\ell\) and the integration parameter as \(s\in[\ell/2,\ell]\).
The fixed lower Ricci bound and fixed diameter \(D_h\) imply
\(\mathcal J_x(\ell,\theta)\leq C(h)\mathcal J_x(s,\theta)\) before the
cut point, by the radial volume-density comparison. This constant is
uniform because \(\ell/s\leq2\) and \(\ell\leq D_h\).
Reverse the \(s,\ell\) integrations; the admissible \(\ell\)-interval has
length at most \(D_h\). Integration over \(x\), and the fact that cut loci
have measure zero, gives \eqref{eq:segment}. All comparison geometry in
this argument concerns \(h\), not the evolving metrics.

Taking \(G=\sqrt{\tr_hg(t)}\), Cauchy--Schwarz and
\eqref{eq:trace-integral} give
\[
 \iint d_{g(t)}(x,y)\,dV_h(x)dV_h(y)
 \leq C_h\left(V\int_X\tr_hg(t)\,dV_h\right)^{1/2}
 =C_h\sqrt m\,V.
\]
Equation \eqref{eq:density}, and \(\Omega=e^{f_0}dV_h\), imply
\(dV_{g(t)}\leq e^{2F}dV_h\).
Transfer both endpoint measures using this inequality and divide by \(V^2\).
The lemma follows, for example, with \(L_0=C_h\sqrt m\,e^{4F}/V\).
\end{proof}

\section{From mean distance to diameter}\label{sec:diameter}

\subsection{Relative volume from pointed Nash entropy}
We recall the relative-volume argument of Jian--Song
\cite[Proposition 3.1, p.~1344]{JS}. Its crucial feature is that it compares
two balls without first bounding their volumes from below in absolute
terms. We retain the coefficient supplied directly by Bamler's gradient
estimate.

We spell out the Riemannian input, with real dimension denoted by \(m\).
Let \(B_t(x,r)\) denote the open ball for \(g(t)\), and let \(\Vol_t\)
denote its volume measure. For a smooth compact Ricci flow, write
\(\N_{x,t}(r^2)\) for pointed
Nash entropy: if the conjugate heat kernel based at \((x,t)\), at time
\(t-r^2\), is
\((4\pi r^2)^{-m/2}e^{-u}\), set
\[
 d\nu=(4\pi r^2)^{-m/2}e^{-u}\,dV_{g(t-r^2)},\qquad
 \N_{x,t}(r^2)=\int_Xu\,d\nu-m/2.
\]
\Needspace{7\baselineskip}
Bamler proves the following in \cite{Bamler},
\begin{itemize}
\item Theorem 5.9, p.~19: for fixed \(s<t\), if \(R(\cdot,s)\geq R_{\min}\),
\[
 \left|\nabla^{g(t)}\N_{\cdot,t}(t-s)\right|
 \leq\sqrt{\frac{m}{2(t-s)}-R_{\min}}.
\]
\item Theorem 6.1, p.~22: if \([t-r^2,t]\) is in the flow interval and
\(R(\cdot,t)\leq r^{-2}\) on \(B_t(x,r)\), then
\[
 \Vol_tB_t(x,r)\geq c_m e^{\N_{x,t}(r^2)}r^m.
\]
\item Theorem 8.1, p.~34: if \(R\geq R_{\min}\) on the full preceding
\(r^2\)-slab, then
\[
 \Vol_tB_t(x,r)\leq U_m(R_{\min}r^2)e^{\N_{x,t}(r^2)}r^m.
\]
The factor corresponding to radius multiplier \(A=1\) has been absorbed
into \(U_m\).
\end{itemize}
None of these three statements assumes a lower entropy bound, a
noncollapsing constant, or a bound for the full curvature tensor.

\begin{lemma}[Jian--Song, from Bamler]\label{lem:relative}
Let \(r>0\). Suppose \([t-r^2,t]\) lies in the flow interval and \(|R|\leq r^{-2}\)
throughout that slab. Then, for any \(y,z\in X\),
\[
 \Vol_tB_t(z,r)\geq
 \kappa_m\exp\left(-a_m\frac{d_t(y,z)}r\right)\Vol_tB_t(y,r),
 \qquad a_m=\sqrt{m/2+1},
\]
where \(\kappa_m>0\) depends only on \(m\).
\end{lemma}
\begin{proof}
Set \(s=t-r^2\) and \(R_{\min}=-r^{-2}\) in the three estimates.
The gradient bound, integrated along a minimizing final-time geodesic,
gives
\[
 \N_{z,t}(r^2)\geq\N_{y,t}(r^2)-a_m d_t(y,z)/r.
\]
Use the upper volume bound at \(y\) and the lower bound at \(z\).
The entropy factors cancel, with \(\kappa_m=c_m/U_m(-1)\).
\end{proof}

\subsection{A first-moment argument}
The following argument adapts the volume-concentration strategy of
Jian--Song \cite[Lemma 4.4 and Proposition 4.1, pp.~1350--1351]{JS}.
Here the first-moment bound supplies a ball of large normalized volume,
and scalar-curvature decay permits comparison at a sufficiently large
fixed radius.

\begin{proposition}\label{prop:diameter}
The diameters of \(g(t)\) are uniformly bounded.
\end{proposition}
\begin{proof}
By averaging the first moment over its center in Lemma~\ref{lem:mean},
at each time there is a point \(y=y_t\) with
\[
 \int_X d_t(y,x)\,d\mu_t(x)\leq L_0.
\]
For \(r>2L_0\), Markov's inequality gives
\(\mu_t(B_t(y,r))>1/2\).
Fix
\[
 \beta_m=\min\{1/4,\tfrac12\kappa_m e^{-3a_m}\}>0,\qquad
 r>\max\{2L_0,L_0/(2\beta_m)\}.
\]
The radius is fixed once and for all. By \eqref{eq:Chow}, for every
sufficiently large \(t\) the interval \([t-r^2,t]\) lies in \([0,\infty)\)
and has \(|R|\leq r^{-2}\).

If some point has distance at least \(3r\) from \(y\), choose a point \(z\)
at distance exactly \(3r\) on a minimizing geodesic. Lemma~\ref{lem:relative}
implies
\[
 \mu_t(B_t(z,r))\geq\beta_m.
\]
Every point in this open ball is farther than \(2r\) from \(y\).
Consequently
\[
 L_0\geq\int_X d_t(y,x)\,d\mu_t(x)
 \geq 2r\,\mu_t(B_t(z,r))
 \geq2r\beta_m>L_0,
\]
a contradiction. Thus every point lies within \(3r\) of \(y\), and
\(\diam(X,g(t))\leq6r\) at all sufficiently late times.
The remaining compact time interval has bounded diameter by smoothness.
\end{proof}

\section{Heat contraction and smooth convergence}\label{sec:convergence}

\subsection{Exponential decay of the Ricci potential}

Let \(D>0\) bound the diameters. Restart the heat equation for \(f\) at
any time \(s\geq0\), and put
\[
 a(s)=\osc_X f(s),\qquad v(x,t)=f(x,t)-\min_Xf(s),\quad t\geq s.
\]
The maximum principle gives \(0\leq v\leq a(s)\).
Using \eqref{eq:bochner} and \(Lv^2=-2q\), we obtain
\[
 L\bigl((t-s)q+v^2/2\bigr)=-2(t-s)|H|^2\leq0.
\]
At time \(s\) this expression is at most \(a(s)^2/2\). Therefore, for \(t>s\),
\[
 |\nabla f(x,t)|\leq\frac{a(s)}{\sqrt{2(t-s)}},\qquad
 a(t)\leq\frac{D}{\sqrt{2(t-s)}}a(s).
\]
Choose the fixed lag \(\ell=2D^2\). Then
\(a(s+\ell)\leq a(s)/2\). Since \(a(t)\) is nonincreasing, iteration and
\eqref{eq:normalization} yield
\begin{equation}\label{eq:exp}
 \|f(t)\|_\infty\leq a(t)\leq C e^{-ct}.
\end{equation}
\subsection{Uniform equivalence and convergence}
In particular
\[
 \sup_{x,t}|\Psi(x,t)|
 \leq2\int_0^\infty\|f(t)\|_\infty\,dt<\infty.
\]
The tensor trace estimate \eqref{eq:Q} now gives a time-independent bound
on \(S\), and \eqref{eq:density} gives
\begin{equation}\label{eq:equivalence}
 N^{-1}\omega_0\leq\omega(t)\leq N\omega_0\qquad(t\geq0).
\end{equation}
Apply Sherman--Weinkove's tensor estimates on fixed coordinate balls and
each unit time slab \([t-1,t]\), \(t\geq1\).
All their inputs, including the slab length and the metric equivalence
constant, are fixed. Hence all derivatives of \(g(t)\) in the \(g_0\)
background are uniformly bounded for \(t\geq1\).

The identity \(f=\log(\Omega/dV_g)\) and \eqref{eq:equivalence} give uniform
bounds for all fixed-background derivatives of \(f\).
On the fixed compact background, the interpolation inequality
\[
 \|u\|_{C^j(g_0)}\leq C_j
 \|u\|_{C^0}^{1/2}\|u\|_{C^{2j}(g_0)}^{1/2}\qquad(j\geq1)
\]
therefore upgrades \eqref{eq:exp} to decay of every derivative.
Taking \(j=k+2\) gives constants \(C_k,c_k>0\) such that
\[
 \|f(t)\|_{C^{k+2}(g_0)}\leq C_k e^{-c_kt}.
\]
As \(\partial_t\omega=-2\ddc f\), the speed is integrable in every \(C^k\)
norm. The forms therefore converge smoothly on the original fixed
manifold, without passing to diffeomorphisms, to a positive K\"ahler form
\(\omega_\infty\). It has class \([\omega_0]\) and
\(\rho(\omega_\infty)=\lim_{t\to\infty}\ddc f(t)=0\).
More quantitatively, integrating the bound for \(\partial_t\omega\) from
\(t\) to infinity gives the exponential estimate in Theorem~\ref{thm:main}
for \(t\geq1\). Enlarging \(C_k\) covers the compact interval \([0,1]\).
Equivalently, \(\Psi\) converges smoothly to a function \(\Psi_\infty\), and
\(\omega_\infty=\omega_0-\ddc\Psi_\infty\). This also makes the equality
of K\"ahler classes explicit.

\section{Uniqueness and the scope of the proof}

If \(\alpha,\beta\) are cohomologous Ricci-flat K\"ahler forms, then
\(\ddc\log(\beta^n/\alpha^n)=0\). Taking the \(\alpha\)-Laplacian and using
compact connectedness shows that the ratio is constant. Equality of total
volumes makes it one. For the positive Hermitian endomorphism
\(A=\alpha^{-1}\beta\),
\[
 \tr A-n-\log\det A=\sum_{j=1}^n(\lambda_j-1-\log\lambda_j)\geq0,
\]
with equality exactly when \(A=I\).
The determinant term vanishes, and cohomology gives
\[
 \int_X\tr_\alpha\beta\,\alpha^n
 =n\int_X\beta\wedge\alpha^{n-1}=n\int_X\alpha^n.
\]
The nonnegative continuous integrand has zero integral, so \(A=I\)
pointwise. Hence \(\alpha=\beta\), completing Theorem~\ref{thm:main}.

The negative of \(\Psi\) is a K\"ahler potential, whose bound here follows
from the diameter estimate and heat contraction. The Schwarz calculation
and tensor derivative estimates retain the classical methods of K\"ahler
geometry. The proof obtains the required global control without invoking
complex Monge--Amp\`ere solvability or global estimates for its scalar
potential equation.

\appendix
\section{Constructing the initial Ricci potential by a linear equation}
\label{app:potential}

For completeness, we recall the part of the compact K\"ahler
\(\partial\bar\partial\)-lemma needed in the proof. It shows explicitly
that constructing \(\Omega\) does not require the Calabi--Yau theorem.

\begin{lemma}\label{lem:ddbar}
Let \(\eta\) be a smooth exact real \((1,1)\)-form on a connected compact
K\"ahler manifold \((X,\omega_0)\) of complex dimension \(n\geq1\). Then
\(\eta=\ddc u\) for a smooth real function \(u\).
\end{lemma}
\begin{proof}
Write \(\Lambda\eta=\tr_{\omega_0}\eta\). Stokes' theorem gives
\[
 \int_X\Lambda\eta\,dV_{g_0}
 =\frac1{(n-1)!}\int_X\eta\wedge\omega_0^{n-1}=0.
\]
The solvability criterion for the scalar Poisson equation on a connected
compact Riemannian manifold therefore gives a smooth real solution of
\(\Delta_{g_0}u=2\Lambda\eta\). Since
\(2\Lambda(\ddc u)=\Delta_{g_0}u\), the real \((1,1)\)-form
\(\sigma=\eta-\ddc u\) is exact, closed, and primitive.
When \(n=1\), primitivity already implies \(\sigma=0\). When \(n\geq2\),
the Hodge star of \(g_0\), with the complex orientation, satisfies
\[
 *\sigma=-\frac{\sigma\wedge\omega_0^{n-2}}{(n-2)!}.
\]
Thus \(d(*\sigma)=0\). Write \(\sigma=d\theta\); another application of
Stokes' theorem yields
\[
 \int_X|\sigma|^2\,dV_{g_0}
 =\int_Xd\theta\wedge *\sigma=0.
\]
Hence \(\sigma=0\), as required.
\end{proof}

Since \([\rho(\omega_0)]=2\pi c_1(X)\) in real cohomology, the hypothesis
of Theorem~\ref{thm:main} makes \(\rho(\omega_0)\) exact. Apply
Lemma~\ref{lem:ddbar} to obtain \(f_0\). The normalization
\(\int_Xe^{f_0}dV_{g_0}=\int_XdV_{g_0}\) then requires only adding a
constant to \(f_0\).\medskip 

\textbf{Acknowledgements.}
The core mathematical results and the final exposition of this paper were generated and refined via iterative prompting of OpenAI models. B.C. thanks OpenAI for access through the ChatGPT for Academic Researchers program. The authors have independently checked all proofs for mathematical correctness.

\end{document}